\documentclass{article}
\usepackage{amsmath,amssymb,amscd,latexsym}
\usepackage{fullpage}
\begin{document}
\parskip 1ex
\parindent0ex
\newtheorem{defn}{Definition}
\newtheorem{lmma}[defn]{Lemma}
\newtheorem{thm}[defn]{Theorem}
\newtheorem{claim}[defn]{Claim}
\newtheorem{cor}[defn]{Corollary}
\newenvironment{proc} {\textsl{Procedure\\}} {\hfill\vskip2ex}
\newtheorem{Thm}[defn]{THEOREM}
\newenvironment{prf}{\noindent{\textsl{Proof\\ }}}{\hfill \rule{1ex}{1.5ex} \vskip2ex}
\newcommand{\nCL}{\rm{CL}}
\newcommand{\nCLo}{\nCL^0}
\newcommand{\ndiv}{\rm{div}}
\newcommand{\Div}{\rm{Div}}
\newcommand{\comment}[1]{}
\title{A new proof of the Nonrationality of Cubic Threefolds}
\author{Tawanda Gwena}
\date{}
\maketitle
\setcounter{section}{-1}
\begin{abstract}
A new proof of the non-rationality of a generic cubic threefold is
given as follows: If a generic cubic threefold were rational then the
associated intermediate Jacobian would be a product of Jacobians of
curves. We degenerate a generic cubic threefold to the Segre Cubic
Threefold and so there is a degeneration of intermediate Jacobians as
well. Associated to the degenerating family of Pryms is a unimodular
system of vectors. Rationality of the generic cubic threefold would
imply that the unimodular system would be cographic dicing. However,
we show that the unimodular system obtained is a well known symmetric
non-cographic dicing called $E_5$.
\end{abstract}

\section{Introduction and History}
The aim of this paper is to use methods of degenerations of Prym varieties
of \cite{ABH} to prove the following theorem: 
\begin{thm}\label{thm1} 
A generic cubic threefold is not birational to $\mathbb {P}^3$.
\end{thm}

There have been several proofs of this theorem. The first proof was by
Clemens and Griffiths in \cite{CG}.  However, even though their proof
does use degenerations, it does not use them directly in the proof of
Theorem \ref{thm1}. So we cannot compare their proof to the others
mentioned here. The advantage which their proof has over the ones
listed here is that it applies to any cubic threefold.

Collino gives a degeneration proof in \cite{CO}. He looks at a family
${\cal X}/S$ of cubic threefolds where $S$ is a smooth but not
necessarily complete curve. The family has the property that for each
$s\in S$ and $ s\ne s_0$ the threefold $X_s$ is smooth, and $X_{s_0}$
has exactly one ordinary double point. To this family he associates a
family of generalized Prym varieties ${\cal A}/S$ such that when $s\ne
s_0$, $A_s=J(X_s)$, the intermediate Jacobian of $X_s$ and $A_{s_0}$
is an extension of a Jacobian $B$ of a curve $C$ by a torus $\tau$.
On the family ${\cal A}/S$ he constructs a relative cartier divisor
$D$ such that when $s\ne s_0$, $D_s$ induces the same polarization as
the theta divisor of $J(X_s)$ and $D_{s_0}=2\Theta$ in $NS(B)$, the
Neron-Severi group. Here $\Theta$ is a divisor which induces a
principal polarization on $B$.  The irrationality is proved by showing
that if all $X_s$ ($s\ne s_0$) were rational then the pair
$(A_{s_0},D_{s_0})$ should be a polarized generalized Jacobian of a
curve with ordinary double points.  The curve $C$ would have to be
both hyperelliptic and of genus 4 embedded in $\mathbb{P}^3$ which is
a contradiction.  This proof works over all algebraically closed
fields where the characteristic is not 2.

Bardelli also does a degeneration proof in \cite{BA}. In his case he
starts out with a family ${\cal X}/ \mathbb {P}^1$ of cubic
threefolds. For $t\ne0$ the fiber $X_t$ is a smooth cubic threefold
and $X_0=\oplus_{j=1}^3X_0^i$ where each $X_0^i$ is isomorphic to
$\mathbb{ P}^3$.  The generalized intermediate jacobian of $X_0$ is an
extension of $\oplus_{j=1}^3J(X_0^i)$ by a torus $\tau$. Then for
$t\ne0$ $H^3(X_t)$ is polarized by a cup product which is denoted
$\Theta(t)$. The family $\{\Theta(t)\}$ specializes to the natural
polarization on $\oplus_{j=1}^3J(X_0^i)$. The Hodge structure
corresponding to $\tau$ is polarized by a bilinear symmetric form
$\psi_{X_0}$. It is shown that $\psi_{X_0}$ is defined in terms of
$\Theta_t$ and the local monodromy of the family of threefolds. By
assuming $X_t$ is rational for $t\ne0$ he shows that $\psi_{X_0}$ is
not the natural polarization on the space of transverse 1-cycles of
any semistable curve. Because of the use of Hodge theory, this proof
only works over $\mathbb{C}$.

The proof in this paper is different in that it is maximal in the
following sense: The Prym variety arising from double covers of stable
curves is an extension of an abelian part by a torus.  For the proofs
above there is an extension of a \textbf{Jacobian} by a torus
$\tau$. In Bardelli's proof $\tau$ is $(\mathbb{C}^*)^2$ (\cite{BA}
6.2.1) and in Collino it is $({k^*})^1$ where $k$ is an algebraically
closed field. For the proof presented here there is no abelian part,
so the torus is $({k^*})^5$ which is the maximum it can be.

The outline of this paper is as follows: In section \ref{pryms} we
gather the results we need from \cite{A} and \cite{ABH} and \cite{V1}.
which we will use later in the paper. Section \ref{uni} presents the
information we need to know about unimodular systems and $E_5$ in
particular. Section \ref{quintics} gives the relationship between
cubic threefolds and Prym varieties. Finally, in section \ref{proof},
we compute the unimodular system for our degeneration, put everything
together and prove Theorem \ref{thm1}.

\textbf{Acknowledgements}: I would like to thank Prof. Valery Alexeev
for his help and attention to my work and Prof. Berndt Sturmfels for
his help with Macaulay2. I would also like to thank Prof. Roy Smith
and Prof. Robert Varley for lots of help and suggestions in the course
of preparing this paper.

\section{Degeneration results}\label{pryms}

In this section we recall the results we need in order for our
degeneration to work. The results are proved in the papers \cite{A},
\cite{ABH} and \cite{N1}. All results work over algebraically closed
fields with characteristic not 2.

\subsection{Jacobians}\label{jac}

Suppose we have a 1-parameter family of smooth curves degenerating to
a stable curve with dual graph $\Gamma$. Then by \cite{A} (or
\cite{ABH}) we get induced data. For us the object of interest is the
cell decomposition obtained (\cite{A} 5.5.) This cell decomposition is
obtained by intersecting the subspace $H_1(\Gamma, \mathbb{R})\subset
C_1(\Gamma,\mathbb{R})$ with the standard cubes in $C_1(\Gamma,
\mathbb{R})$. This decomposition, called a \textbf{cographic dicing},
does not depend on the 1-parameter family and can be obtained from a
unimodular system (\cite{ABH}, 2.3 (J6)).

\subsection{Pryms}\label{prym}

Suppose we have a 1-parameter family of curves $(\tilde{C}, \iota)/S$
such that the generic fiber is a smooth curve with a base-point-free
involution and the degenerate curve $(\tilde{C}_0, \iota)$ is a stable
curve. Then degenerate data can be obtained for this family
(\cite{ABH}, Section 2.4 PP0--PP6). In particular the data we need is
the cell decomposition of PP6 which we now describe.

For the rest of this section we will the drop the subscript $0$ from
$\tilde{C}_0$ and $\tilde{\Gamma}_0$.  Let the dual graph of
$\tilde{C}$ be $\tilde{\Gamma}$.  To get this cell decomposition we do
the following: We define a map
\begin{eqnarray*}
\pi^-:H_1(\tilde{\Gamma},\mathbb{ Z})&\longrightarrow
&H_1(\tilde{\Gamma}, \mbox{$\frac{1}{2}$}\mathbb{ Z})\\ h &\longmapsto
&\mbox{$\frac{1}{2}$}(h-\iota(h))
\end{eqnarray*}
 
Let $X^-:=\pi^-(H_1(\tilde{\Gamma},\mathbb{ Z}))$.  The space
$X^-\otimes \mathbb{R}$ is contained in $C_1(\tilde{\Gamma},
\mathbb{R})$. Each edge $e_j$ of $\tilde{\Gamma}$ defines a coordinate
function $z_j$ in $C_1(\tilde{\Gamma}, \mathbb{R})$. Let $m_j=1$ if
$z_j:X^-\rightarrow \mathbb{Z}$ is surjective and $m_j=2$ if
$z_j:X^-\rightarrow \frac{1}{2}\mathbb{Z}$ is surjective. The
functions $m_jz_j$ define a cell decomposition of $X^-\otimes
\mathbb{R}$.

The cell decomposition defined on $X^-\otimes \mathbb{R}$ is
independent of the 1-parameter family if the vertices of the cell
decomposition are precisely the points of $X^-$, i.e. the linear
functions define a dicing of $X^-$. This is condition (*) in
\cite{ABH}. An equivalent, easier to verify, condition is given in
\cite{V1} as Theorem 0.1. This condition can be stated as follows:

\begin{lmma}[\cite{V1} Theorem 0.1]\label{vit}
The cell decomposition depends only on the degenerate fiber
$(\tilde{C},\iota)$ if the following is true: There do not exist two
connected subgraphs $\Gamma_0, \Gamma_1$ of the dual graph of
$\tilde{C}$ such that $\iota(\Gamma_i)=\Gamma_i$, $(i=0,1)$ and there
are at least four edges connecting $\Gamma_0$ and $\Gamma_1$.
\end{lmma} 

\section{Unimodular systems}\label{uni}
\begin{defn} [\cite{DG}]
A system of $m\ge n$ vectors $R$ generating $\mathbb {R}^n$ is called
a unimodular system ($U$-system) if when we write any $n$ vectors of
$R$ as columns in terms of a basis $B\subset R$ we obtain a matrix
which is totally unimodular, i.e the maximal minors are either $0$,
$1$ or $-1$.
\end{defn}

Suppose $R$ is a set of vectors spanning $\mathbb {R}^n$. Then these
vectors define a family $H(R)$ of of parallel hyperplanes
$H(\textbf{r},z)=\{\textbf{x}\in \mathbb{ R}^n :
\textbf{x}\cdot\textbf{r}=z\}, z\in \mathbb{ Z},\textbf{r}\in R$. If
$B\subset R$ is a basis for $\mathbb{ R}^n$ then the intersection
points of hyperplanes in $H(B)$ is a lattice and $H(R)$ is then called
a {\bf lattice dicing}. The set of intersection points of hyperplanes
of $H(R)$ is a lattice if and only if $R$ is a unimodular system.

There are unimodular systems which are not cographic.
In our case, we are interested in a very nice and exceptional one called
$E_5$. It is represented by the matrix 
$$
\left(\begin{array}{ccccc|ccccc}
1& 0& 0& 0& 0& 1 & 0 & 0 & 1 & 1\\
0& 1& 0& 0& 0& 1 & 1 & 0 & 0 & 1\\
0& 0& 1& 0& 0& 0 & 1 & 1 & 0 & 1\\
0& 0& 0& 1& 0& 0 & 0 & 1 & 1 & 1\\
0& 0& 0& 0& 1& 1 & 1 & 1 & 1 & 1 
\end{array}\right).
 $$ 

The fact that it is not cographic is proved as Cor 13.2.5 in \cite{O1}.

\section{Cubic Threefolds and Plane Quintics}\label{quintics}
In this section we show the connection between cubic threefolds and
plane quintics for nonsingular threefolds and the Segre threefold.
Most of the information here is put together from the following
sources: in \cite{D} Sections 4.8 and 5.17, \cite{H} Section 3.1 and
3.2, \cite{M1}, \cite{M2}, \cite{CG} and \cite{SR}.
\subsection{Smooth Cubic Threefolds}
Let $X$ be a smooth cubic threefold in $\mathbb{ P}^4$. The lines in
$X$ form a surface $F$ called a Fano surface (\cite{M1}, 1.1;
\cite{CG}, Theorem 7.8). Pick a generic enough line $\ell$ in
$F\subset X$ (satisfying conditions in \cite{M1}, Prop 1.25). The
space of planes through $\ell$ is parametrized by a projective space
$Y=\mathbb{ P}^2$. Let $C_\ell\subset Y$ be the planes $L$ such that
$L\cap X$ consists of three lines. Then $C_\ell$ is a plane quintic
curve. Also, let $\tilde{C_\ell}=\{\ell'\in F|\ell\cap\ell'\ne
\emptyset \}$. Then $\tilde{C_\ell}$ is a curve in $F$ and there is a
natural map $q:\tilde{C_\ell}\rightarrow C_\ell$. We get an involution
$\iota$ on $\tilde{C}$ as follows: for a point $L\in C$
$q^{-1}(L)=\{\ell',\ell''\}$. So for $\ell'\in \tilde{C}$ we have
$\iota(\ell')=\ell''$. From this double cover we obtain a Prym variety
$P(\tilde{C_\ell}, \iota)$. However the Prym variety is independent of
$\ell$, so we can write $P(X)$. This construction is in \cite{M2}
Section 2, and \cite{CG}, Appendix C.

\subsection{The Segre Threefold} \label{segre}
Now we construct a singular plane quintic and a double cover.
\begin{defn}
The Segre cubic threefold ${\cal S}$ is given by the following
equation in $\mathbb{P}^5$:
$$ \left\{(x_0,\ldots, x_5):\sum_{j=0}^5 x_j^3=\sum_{j=0}^5
x_j=0\right\}.$$
\end{defn}

This equation only works with fields of characteristic not 3. In all
characteristics not 2 the following works:
$$ \left\{(x_0,\ldots, x_4):\sum_{\mbox{\shortstack{$i$, $j$, $k$\\ distinct}}}2x_i x_j
x_k+\sum_{i\ne j}x_i^2 x_j =0,\ \ \  i,j,k\in \{0,\ldots, 4\}\right\}$$

The Segre threefold is the unique (up to isomorphism) threefold with
ten nodes. The equation above shows that the cubic threefold is
invariant under the action of the symmetric group $S_6$ on the
coordinates. One node is given by the coordinates $(1,1,1,-1,-1,-1)$
and the other nine are obtained by the action of $S_6$ on the
coordinates.

We can obtain ${\cal S}$ as follows: Let $p_1,\ldots, p_5$ be points
in general position in $\mathbb{P}^3$ and for each pair of points
$\ell_{ij}$ is the line joining $p_i$ and $p_j$.  Blow up
$\mathbb{P}^3$ at the five points and then blow down the proper
transforms of the ten lines $\ell_{ij}$. The image is $S$ in
$\mathbb{P}^4$.

There are $15$ planes contained in ${\cal S}$. Ten of them come from
proper transforms of the planes in $\mathbb{P}^3$ containing three
points $p_i$, $ p_j$ and $p_k$. These are labeled $\Pi_{ijk}$ ($i< j<
k$). The other five are images of the exceptional divisors coming from
the blowup. These are labeled $\Pi_{i}$ ($1\le i\le 5$). Each of the
15 planes contains 4 nodes.

${\cal S}$ contains 6 two dimensional families of lines, $R_i$ ($1\le i\le
5$) which are proper transforms of lines through $p_i$ and $R_0$ which
is the family of twisted cubics through the five points. Each line in
$R_i$ goes through the five planes $\Pi_{ijk}$ ($i < j < k$). Each line in
$R_0$ goes through $\Pi_i$ (\cite{SR} VIII 2.32).

We choose a line $\ell$ in $R_1$ which 
does not go through a node and goes through exactly 5
planes. We project from
$\ell$ onto $\mathbb{P}^2$. 
The images of each of the five planes is in the degeneracy locus for
${\cal S}$. Each line meets each of the other four lines at a node. Each node is where
the preimage is a plane containing a node on ${\cal S}$.
The five lines form a pentagon in $\mathbb{P}^2$ which we call
$C_0$. The dual graph of $C_0$ is the complete graph on five vertices $K_5$. 

As in the smooth case, we can construct a double cover for the curve $C$ above.
This double cover will lie in $G(2,5)$. The preimage of each point in
$C$ consists of two lines, excluding $\ell$, which are points in
the double cover. The whole cover can be described as follows:
Take ten copies of $\mathbb{ P}^1$, $L_i^\epsilon$ where $1\le i\le 5$
and $\epsilon$ is either 0 or 1. Each line $L_i^\epsilon$ has four
points marked on it $p^\epsilon_{i,j}$ where $1\le j\le 5$ and $j\ne
i$.  In the following the notation $(p^0_{i,j} \sim p^1_{j,i})$ means
that the points $p^0_{i,j}$ and $p^1_{j,i}$ are identified.
The double cover $\tilde{C}$ is  
$$
\tilde{C} = \left.\left(\coprod_{i,\epsilon} L_i^\epsilon\right) \right/ \left(p^0_{i,j} \sim p^1_{j,i} \right)$$

The dual graph  $\Gamma_0=\Gamma(\tilde{C_0})$ is shown below. If the vertices of
$\Gamma(C_0)$ are labelled $v_j$, the the vertices $a_j$ and $b_j$ both
map to $v_j$. The map of edges is given by the map of vertices.\\

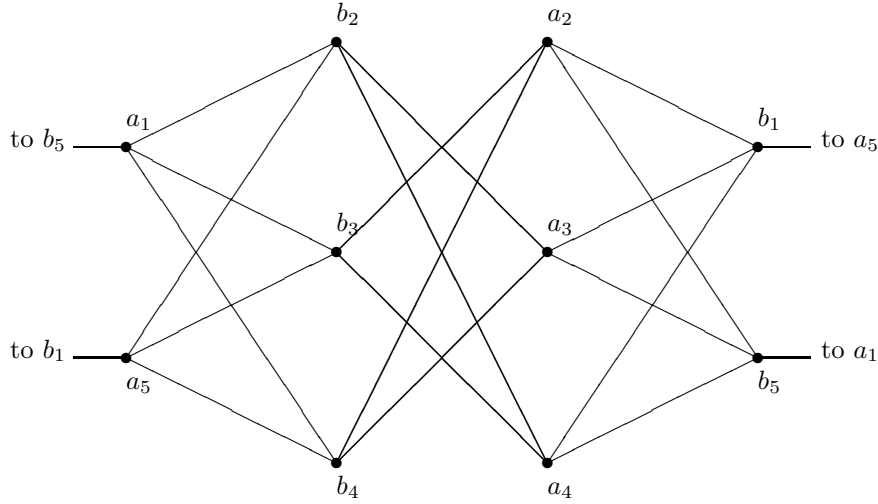
\begin{center}
\setlength{\unitlength}{0.7mm}
\begin{picture}(200,90)(-100,-50)
\put(-60,20){\line(2,1){40}}
\put(-60,20){\line(2,-1){40}}
\put(-60,20){\line(2,-3){40}}

\put(-60,-20){\line(2,-1){40}}
\put(-60,-20){\line(2,1){40}}
\put(-60,-20){\line(2,3){40}}

\put(-20,40){\line(1,-1){40}}
\put(-20,40){\line(1,-2){40}}
\put(-20,-40){\line(1,1){40}}
\put(-20,-40){\line(1,2){40}}
\put(-20,0){\line(1,-1){40}}
\put(-20,0){\line(1,1){40}}

\put(60,20){\line(-2,1){40}}
\put(60,20){\line(-2,-1){40}}
\put(60,20){\line(-2,-3){40}}

\put(60,-20){\line(-2,-1){40}}
\put(60,-20){\line(-2,1){40}}
\put(60,-20){\line(-2,3){40}}

\put(20,40){\line(-1,-1){40}}
\put(20,40){\line(-1,-2){40}}
\put(20,-40){\line(-1,1){40}}
\put(20,-40){\line(-1,2){40}}
\put(20,0){\line(-1,-1){40}}
\put(20,0){\line(-1,1){40}}
\put(-70,20){\line(1,0){10}} \put(-82,20){to $b_5$}
\put(-70,-20){\line(1,0){10}}\put(-82,-20){to $b_1$}
\put(70,20){\line(-1,0){10}}\put(72,20){to $a_5$}
\put(70,-20){\line(-1,0){10}}\put(72,-20){to $a_1$}

\put(-60,24){$a_1$} \put(-60,20){\circle*{2}}
\put(-60,-26){$a_5$}\put(-60,-20){\circle*{2}}
\put(-20,44){$b_2$}\put(-20,40){\circle*{2}}
\put(-20,4){$b_3$}\put(-20,0){\circle*{2}}
\put(-20,-46){$b_4$}\put(-20,-40){\circle*{2}}

\put(60,24){$b_1$}\put(60,20){\circle*{2}}
\put(60,-26){$b_5$}\put(60,-20){\circle*{2}}
\put(20,44){$a_2$}\put(20,40){\circle*{2}}
\put(20,4){$a_3$}\put(20,0){\circle*{2}}
\put(20,-46){$a_4$}\put(20,-40){\circle*{2}}
\end{picture}\\

\begin{minipage}[c]{80ex}
\small
The edges are named as follows: 
$e_1=(b_3,a_2)$, $e_2=(a_4,b_2)$, $e_3=(a_5,b_3)$,  
$e_4=(a_5,b_4)$, $e_5=(a_5,b_1)$, $e_6=(a_4,b_3)$,
$e_7=(b_3,a_1)$, $e_8=(b_2,a_5)$, $e_9=(a_1,b_4)$, $e_{10}=(b_2,a_1)$. 
The rest of the edges are named as follows: if edge $e_i$ is
$(a_j,b_k)$ ($(b_j,a_k)$ respectively) the edge $e_i'$ is $(b_j,a_k)$,
($(a_j,b_k)$ respectively).

The tree used to form the basis of $H_1(\Gamma, \mathbb{Z})$ is given by the edges
 $e_6=(a_4,b_3)$,  $e_7=(b_3,a_1)$, $e_8=(b_2,a_5)$,
$e_9=(a_1,b_4)$, $e_{10}=(b_2,a_1)$ and
  $e_{10}'=(a_2,b_1)$, $e_{7}'(a_3,b_1)$,
$e_{9}'=(b_1,a_4)$, $e_{8}'=(a_2,b_5)$.

\end{minipage}
\end{center}

\section{Proof of the Main Theorem}\label{proof}

We now prove the irrationality of cubic threefolds by making
use of the following lemma which relates cubic threefolds and Prym varieties.

\begin{lmma}[\cite{M2} Thm 3.11]\label{mur1}
Let $\mbox{char}(k)\ne 2$. Let $X$ be a nonsingular cubic threefold in
$\mathbb{ P}^4$, defined over $k$. If there exists a birational
transformation between $X$ and $\mathbb{ P}^3$ then the canonically
polarized prym variety $(P(X), \Xi)$ associated with $X$ is
isomorphic, as a polarized abelian variety, to a product of
canonically polarized Jacobian varieties of curves.
\end{lmma}

Given a smooth plane quintic curve with a double cover we degenerate
it to the stable curve $C_0$ above. The generalized Prym
$P(\tilde{C},\iota)$ can be easily shown to be isomorphic to
$({k^*})^5$.

Using the results on degeneration above we can compute the unimodular
system for the Delaunay decomposition of the degeneration.  By
constructing the dual graph of the double cover and using the
algorithm in \cite{ABH}, outlined in Section \ref{prym} we prove the
following theorem:

\begin{thm}\label{42}
The unimodular system for cell decomposition associated to a family of
cubic threefolds degenerating to the Segre Threefold is $E_5$.
\end{thm}

Before we do the proof we need to verify that if a family of cubic
threefolds degerate to ${\cal S}$ the the family of double covers satisfies
Lemma \ref{vit}.

Let $\Gamma$ be as above. Suppose we have the two subgraphs $\Gamma_1$
and $\Gamma_2$ of $\Gamma$.  Then the first possible case is that
$\Gamma_1$ has 2 vertices and $\Gamma_2$ has 8 vertices. This is not
possible because then $\Gamma_1$ would not be connected. Suppose the
vertices in $\Gamma_1$ correspond to the lines $L_1^0$ and $L_1^1$
(using the notation from above).  The line $L_1^0$ and $L_1^1$ do not
meet, so on the dual graph their corresponding vertices do not have an
edge between them.  The second possible case is if $\Gamma_1$ and
$\Gamma_2$ have 4 and 6 vertices respectively.  Suppose without loss
of generality the lines $L_1^0$,$L_1^1$, $L_2^0$ and $L_2^1$ are in
$\Gamma_1$. This would imply $\Gamma_1$ is not connected because the
connected subgraph with $L_1^0$ and $L_2^1$ is not connected with the
subgraph with $L_1^1$ and $L_2^0$. So the Delaunay decomposition is
independent of the 1-parameter family, therefore we can obtain a
unimodular system.

\begin{prf}
[of Theorem \ref{42}]

We now compute $X^-\subset C_1(\Gamma, \mathbb{Z})$.

The following is a basis for $H_1(\Gamma,\mathbb{ Z})$.

\begin{minipage}{40ex}
\begin{eqnarray*}
h_1&=&e_6'+e_7'-e_9'+e_6+e_7-e_9\\
h_2&=&e_1'+e_7'+e_9'+e_6+e_7+e_{10}\\
h_3&=&e_1-e_{10}'-e_9'+e_6\\
h_4&=&e_2+e_6+e_7+e_{10}\\
h_5&=&e_2'-e_{10}'-e_9'+e_6+e_7-e_9\\
\phantom{.}
\end{eqnarray*}
\end{minipage}
\begin{minipage}{40ex}
\begin{eqnarray*}
h_6&=&e_5'-e_8'-e_{10}'-e_9'+e_6+e_7\\
h_7&=&e_5-e_9'+e_6+e_7+e_{10}+e_8\\
h_8&=&e_4+e_9+e_{10}-e_8\\
h_9&=&e_3+e_7+e_{10}-e_8\\
h_8'&=&e_4'+e_9'+e_{10}'-e_8'\\
h_9'&=&e_3'+e_7'+e_{10}'-e_8'.
\end{eqnarray*}
\end{minipage}

The basis for $X^-$ is as follows

\begin{eqnarray*}
\ell_1&=&\frac{1}{2}(h_2-\iota(h_2)) = \frac{1}{2}(h_3-\iota(h_3))\\
\ell_2&=&\frac{1}{2}(h_4-\iota(h_4)) = \frac{1}{2}(h_5-\iota(h_5))\\
\ell_3&=&\frac{1}{2}(h_9-\iota(h_9)) = \frac{1}{2}(h_9'-\iota(h_9'))\\
\ell_4&=&\frac{1}{2}(h_8-\iota(h_8)) = \frac{1}{2}(h_8'-\iota(h_8'))\\
\ell_5&=&\frac{1}{2}(h_6-\iota(h_6)) = \frac{1}{2}(h_7-\iota(h_7))
\end{eqnarray*}

and $X^-=<\ell_1,..., \ell_5>$. The unimodular system is obtained
by seeing how the edges $e_j$ restrict to $X^-$. The unimodular matrix for the dicing
of $X^-\otimes \mathbb{R}$ is $(a_{ij})$
where $a_{ij}$ is defined to be $1$ if $\ell_i$ contains $e_j$. 
This matrix is $E_5$.
\end{prf}

We are now in a position to prove Theorem \ref{thm1}.

Let $\ell_0$ be a line in ${\cal S}$ (chosen as in Section
\ref{segre}). Choose a smooth cubic threefold $X$ such that it also
contains $\ell_0$. We have a pencil of cubic threefolds which
contain $\ell_0$: $$X_{a,b}=aX+b{\cal S}\ \ \ \ \ \ \ \ (a:b)\in
\mathbb{P}^1.$$

By restricting to some open subset $S$ of $\mathbb{P}^1$ we get a family
${\cal X}/S$ of cubic threefolds. For each $s$ there is a Prym
$P(X_s)$. If the threefolds in our family were rational then by
Section \ref{jac} and Lemma \ref{mur1} the family of Pryms we obtain
should give a cographic unimodular system. 
But by Theorem \ref{42} we get $E_5$ which we know
 is not cographic. So our original supposition that the
cubic threefolds were rational is false. \hfill \rule{1ex}{1.5ex}

\end{document}